\documentclass[12pt]{article}
\usepackage[dvips]{graphicx}
\usepackage{latexsym}
\usepackage{amssymb}

\newcommand{\NN}{\mathbb{N}}
\newcommand{\C}{\mathbb{C}}
\newcommand{\CP}{\mathbb{CP}}
\newcommand{\RP}{\mathbb{RP}}
\newcommand{\R}{\mathbb{R}}

\renewcommand{\d}{\mathrm{d}}

\newcommand{\koniec}{\begin{flushright}  $\Box $ \end{flushright}}
\def\be{\begin{equation}}
\def\ee{\end{equation}}
\def\theequation{\thesection.\arabic{equation}}

\def\Sm{\Sigma}

\def\Om{\Omega}

\def\G{\Gamma}
\def\OO{{\cal O}}
\def\om{\omega}

\def\p{\partial}

\def\sm{\sigma}

\newcommand{\hook}{{\setlength{\unitlength}{11pt}   
                   \begin{picture}(.833,.8)
                   \put(.15,.08){\line(1,0){.35}}
                   \put(.5,.08){\line(0,1){.5}}
                   \end{picture}}}
\def\a{\alpha}
\def\ll{{\lambda}}

\newtheorem{theo}{Theorem}[section] 
\newtheorem{prop}[theo]{Proposition}

\newtheorem{col}[theo]{Corollary}

\begin{document}
\title{Harmonic functions, central quadrics, and twistor theory}
\author{Maciej Dunajski\\
Department of Applied Mathematics and Theoretical Physics, \\
Cambridge University,\\
Wilberforce Road, Cambridge, CB3 OWA, UK
}
\date{} 
\maketitle
\abstract {
Solutions to the $n$-dimensional Laplace equation 
which are constant on a central quadric are found.
The associated twistor description of the
case $n=3$ is used to characterise Gibbons-Hawking metrics 
with tri-holomorphic $SL(2, \C)$ symmetry. 
} 
\noindent
\section{Introduction}
Let $V(x_1, ..., x_n)\in\C$ be a solution to a PDE
\be
\label{PDE0}
F(V, V_i, V_{ij}, ..., V_{ij\cdots k}, x_i)=0,
\ee
where $V_i=\p V/\p x_i$, and $(x_1, ..., x_n)\in\C^n$. 
The usual way of reducing this PDE to an ODE 
is to determine a group of transformations acting on dependent and independent
variables, such that $V$ is transformed to a different 
solution of (\ref{PDE0}), and reduce the number of independent variables
down to one.  There are several algorithmic procedures of 
varying levels of sophistication
(Noether's theorem, Lie Point symmetries, Cartan-K{\"a}hler formalism, ...)
to study  such reductions. See \cite{O95} for a very good review of some of
these methods.

It is interesting to seek non-symmetric ways of reducing PDEs to ODE.
One such method is the hyper-surface ansatz. Let $\Sm\subset\C^n$ be an 
algebraic  hyper-surface. The ansatz 
(which is motivated by the work of Darboux \cite{Darboux}
on orthogonal curvilinear coordinates) 
is to seek solutions constant on $\Sm$, 
or equivalently to seek symmetric objects 
\[
M(V), M^i(V), M^{ij}(V), ..., M^{ij...k}(V), 
\] 
so that
a solution of equation (\ref{PDE0})  is determined implicitly 
by
\be
\label{qq2}
Q(x_i, V):=M(V) + M^{i}(V)x_i + M^{ij}(V)x_ix_j+ ...+M^{ij...k}(V)x_ix_j...x_k =C,
\ee
where $C$ is a constant. Here $\Sm$ should be regarded as the zero locus of a 
polynomial $Q(x_i, V)-C$ in $\C^n$.

If $V$ satisfies (\ref{PDE0}) and  the algebraic 
constraint (\ref{qq2}), then so does $g^t(V)$, there $g^t$ is a flow
generated by any section of $T\Sm$. 
Note however that vectors tangent to $\Sm$ do not generate symmetries
of (\ref{PDE0}), as the choice of $\Sm$ depends on $V$.

In this paper we shall look the quadric ansatz
\be
\label{q2}
Q(x_i, V):=M^{ij}(V)x_ix_j=C,
\ee
which is the simplest non-trivial case of (\ref{qq2}).
The ansatz  can be    made whenever we have a PDE of the form
\be
\label{PDE}
\frac{\p}{\p x_j}\Big( \eta_{ij}(V)\frac{\p V}{\p x_i}\Big)=0,
\ee
where  $\eta$ is a given symmetric matrix whose components depend on $V$.

The quadric ansatz
has been applied to two non-linear integrable PDEs:
If $V$ is a solution of $SU(\infty)$ Toda equation
\[
\exp{(V)}V_{12}-V_{33}=0
\]
then $M$
can be determined in terms of the third Painlev\'e transcendent
\cite{T95}. If $V$ is a solution to dispersionless KP equation 
\[
(V_1-VV_2)_2-V_{33}=0
\]
then
$M$  can be  determined by one of the first two Painlev\'e transcendents 
\cite{DT02}.
Both $SU(\infty)$ Toda, and dKP equations are solvable by twistor transform,
and it would be desirable to achieve a 
characterisation twistor spaces corresponding to  
solutions constant on central quadrics.

Motivated by this problem 
we shall apply the quadric ansatz to the Laplace equation.
We shall work in the holomorphic category, and regard $V$ as a holomorphic
function of complex coordinates $x_1, ..., x_n$. We shall abuse the terminology, and
call $V$ a harmonic function, whenever it satisfies a complexified Laplace 
equation.
The reality conditions may be imposed if desired, to characterise 
real solutions in signatures $(n-r, r)$ when $r=1, ..., n$.

In the next section we shall find harmonic functions constant
on central quadrics (Theorem \ref{main_theo}).
In the remaining sections we shall focus on the three--dimensional case.
In section 3  harmonic functions will be related to solutions
to SDiff$(\Sm)$ Nahm's equations  by means of a 
hodograph transformation. Here SDiff$(\Sm)$ is a group of holomorphic 
symplectomorphisms 
of a two-dimensional complex symplectic manifold $\Sm$.
Harmonic functions constant on a central quadric will be characterised
by a reduction form  SDiff$(\Sm)$ to $SL(2, \C)$ (Proposition \ref{CCONE}).
In section 4 we shall review a twistor construction of 
solutions to $3D$ Laplace equation, and characterise 
solutions 
constant on central quadrics  in terms of $SL(2, \C)$ invariant holomorphic 
line 
bundles over $T\CP^1$ (Theorem \ref{Twistor_th}). 
In section 5 we shall characterise $\C^*$ invariant 
complexified hyper-K\"ahler metrics in $4D$ which admit 
tri-holomorphic $SL(2, \C)$ transitive action 
(Proposition \ref{Gibbons_Hawking}). In section 6 we shall give an example
illustrating the construction.
The basic facts about bundles over $\CP^1$ used in the paper are 
collected in the Appendix A. Appendix B is devoted to SDiff$(\Sm)$
Nahm's equations.
\section{Quadric ansatz for Laplace equation}
\begin{theo}
\label{main_theo}
Solutions of the Laplace equation
\be
\label{laplace}
\sum_{i=1}^n\frac{\p^2 V}{\p x_i^2}=0
\ee
constant on a central
quadric are given by
\be
\label{new_sol}
V=\int\frac{\d H}{\sqrt{(H-\beta_1)(H-\beta_2)\cdots(H-\beta_n)}},
\ee
where 
\be
\label{quadD}
\sum_{i=1}^n\frac{{x_i}^2}{H-\beta_i}=C,
\ee 
and $C, \beta_1, \beta_2, ..., \beta_n$ are constants (which can be normalised 
so  that $\beta_1+\beta_2+...+\beta_n=1$).
\end{theo}
{\bf Proof.}
Equation (\ref{laplace}) is equivalent to (\ref{PDE}) if
$\eta_{ij}$ is an identity matrix. We assume that the level sets of $V$ are of the form  {\em(\ref{q2})}, and we differentiate (\ref{q2}) implicitly to find
\be
\label{ux}
\frac{\p V}{\p x_i}=-\frac{2}{\dot{Q}}M^{ij}x_j,\qquad
\mbox{where}\qquad \dot{Q}=\frac{\p Q}{\p V}.
\ee
Now we substitute this into (\ref{PDE}) and integrate once with
respect to $V$. Introducing $g(V)$ by
\be
\label{gu}
\dot{g}=\frac{1}{2}\eta_{ij}M^{ij}=\frac{1}{2}\mbox{trace}\;{(\bf {\eta}M)}
\ee
we obtain
\[
(g\dot{M}^{ij}-M^{ik}\eta_{km}M^{mj})x_ix_j=0,
\]
so that as a matrix ODE
\be
\label{mode}
g\dot{\bf{M}}={\bf{M{\eta}M}}.
\ee
This equation simplifies if written in terms of another matrix
${\bf N}(V)$ where
\be
\label{NM}
{\bf N}=-{\bf M}^{-1}
\ee
for then
\be
\label{gnb}
g\dot{\bf{N}}={\bf \eta},
\ee
and $g$ can be given in terms of $\triangle=\det{({\bf{N}})}$ by
\be
\label{const}
g^2\triangle=\zeta=\mbox{constant}.
\ee
Equation (\ref{gnb}) implies that ${\bf N}(V)$ can
be written as
\be
\label{NdKP}
{\bf N}(u)=\left (
\begin{array}{ccccc}
H_{11}&h_{12}&h_{13}&...&h_{1n}\\
h_{21}&H_{22}&h_{23}&...&h_{2n} \\
...&...&...&...&...\\
...&...&...&...&...\\
h_{n1}&h_{n2}&...&h_{(n-1)n}&H_{nn}
\end{array}
\right )
\ee
where $h_{ij}$ are constants, while
\be
\label{XYZ}
\dot{H}_{ii}=g^{-1}, \qquad H_{11}-H_{ii}=\gamma_i=\mbox{constant},
\qquad i=1, ..., n.
\ee
The constants $h_{ij}$ can be eliminated by $SO(n, \C)$ transformations
of (\ref{laplace}).
Let $H_{11}=H(V)$. Equation (\ref{const}) yields 
\[
H(H+\gamma_2)...(H+\gamma_n)=\zeta \dot{H}^2,
\]
and ${\bf M}=-$diag$[H^{-1}, (H+\gamma_2)^{-1}, ..., (H+\gamma_n)^{-1}]$
so the quadric (\ref{q2}) is diagonal.
We can rescale $H$, shift it by a constant and define a new set of constants
$\beta_1, ..., \beta_n$
so that (\ref{q2}) yields (\ref{quadD})
where  $H(V)$ satisfies
\be
\label{triangle}
\prod_{i=1}^n(H-\beta_i)=\dot{H}^2.
\ee
The last equation is readily solved giving (\ref{new_sol}).
\koniec
{\bf Remarks}
\begin{enumerate}
\item In general solutions to the Laplace equations obtained from the quadric 
ansatz with $C\neq 0$ 
do not admit any symmetries. Solutions which are constant 
on a central cone ($C=0$) 
are invariant under scaling transformations $x_i\rightarrow sx_i$. 
\item 
If V is a harmonic function  given by the
quadric ansatz (\ref{q2}) then $\hat{V}=\p V/\p C$ is also harmonic 
(but not necessarily constant on a quadric).
Repeating the process yields an infinite 
set \[V, \qquad\p V/\p C,\qquad  \p^2 V/\p C^2,\qquad ...\] of solutions associated to the quadric ansatz. They can all be found by 
implicit differentiation of 
(\ref{q2})
\[
\dot{Q}\frac{\p V}{\p C}=C, \qquad
\ddot{Q}\Big(\frac{\p V}{\p C}\Big)^2+\dot{Q}\frac{\p^2 V}{\p C^2}=0, \qquad ...
\]
where $\dot{Q}={\dot{M}^{ij}x_ix_j}$.
For example 
\be
\label{exp_V}
\hat{V}=
-\Big(\prod_{i=1}^n(H-\beta_i)\Big)^{-1/2}
\Big(\sum_{i=1}^n\frac{{x_i}^2}{(H-\beta_i)^2}\Big)^{-1},
\ee
where $H$ is an algebraic root of (\ref{quadD}), and there is no need for
hyper-elliptic integrals! Implicit differentiation of (\ref{new_sol}), and setting $c=1$ shows that $2\hat{V}=\Upsilon(V)$, where $\Upsilon=x_i(\p/\p x_i)$ 
is the Euler's homogeneity operator.
\item
We can impose the Euclidean reality conditions, and seek solutions constant on 
confocal ellipsoids. 

For $n=2$ the  solutions can be written in terms of holomorphic functions. The solution  (\ref{new_sol}) is $\Re f(z)$, where $z=x+\sqrt{-1}y$, and
$f=\ln(z+\sqrt{z^2+\a})^2$ for $\a=(2\beta_1-1)/C$.

The case $n=3$ of (\ref{exp_V}) has been previously 
characterised \cite{M70} in elliptic coordinates,
and revisited in \cite{G03} in a context of gravitational instantons.
 The description of the 
quadric ansatz in terms of arbitrary holomorphic functions 
is also possible for $n=3$, by means of twistor theory. This will be done in 
the Section 4.
\end{enumerate}
\section{The $SL(2, \C)$ Nahm equation  and the quadric ansatz}
In this section we shall present a  hodograph transformation between 
the Laplace equation and a system of first order PDEs. In three dimensions
(when the first order system  is SDiff$(\Sm)$ Nahm's equations) solutions
constant on quadrics will be characterised by a choice of
$SL(2, \C)\subset$SDiff$(\Sm)$.

Equation (\ref{laplace}) is equivalent to 
\be
\label{gdkp}
\d\ast \d V=0,
\ee 
where the Hodge operator $\ast$ defined by the metric $\d {x_1}^2+...+\d {x_n}^2$, and
the volume form $\d x^1\wedge\cdots\wedge\d x^n$.
We say that $V$ generic if $|\d V|^2=\d V\wedge\ast\d V\neq 0$.
The next Proposition shows that in the generic case the Laplace equation
is equivalent to a system of first order PDEs, which we propose to call
the Nambu-Nahm equations
\begin{prop}
\label{Nambu}
Let $V=V(x_i)$ be a generic solution to the Laplace equation
{\em (\ref{laplace})}. The functions $P_a(x_i),a=1,\cdots n-1$ can be found
such that  
$
x_i=x_i(V, P_a)
$
satisfy the following set of first order PDEs:
\be
\label{new_equ}
\dot{x}_i=\{x_1, ..., x_{i-1},x_{i+1}, ..., x_n\}_{NB}, 
\ee
where the Nambu bracket appearing on the RHS of {\em(\ref{new_equ})} is
given in terms of a Jacobian determinant 
by $f^{-1}\p (x\neq x_i)/\p(P)$, 
and $f(P_a)$ is a nowhere vanishing function which locally can be set to $1$.

Conversely, if $
x_i=x_i(V, P_a) $ solve  {\em(\ref{new_equ})}
then $V(x_i)$  satisfies {\em(\ref{laplace})}.
\end{prop}
{\bf Proof.}
In the generic case $|\d V|\neq 0$ the condition
(\ref{gdkp}), and the Poincare Lemma 
imply the existence of a foliation $U=\C\times
\Sm$ of $U\subset \C^n$ by $(n-1)$-dimensional complex manifolds
$\Sm$ with a holomorphic volume form $\om$ such that
\be
\label{poincare}
\ast\d V=\om.
\ee
Let $P_a=(P_1, ..., P_{n-1})$ be local holomorphic coordinates\footnote{
In the null case $|\d V|=0$ the 
existence of  $P_1, ..., P_{n-1}$ can not be deduced
even locally.} on $\Sm$, and
let 
\[
\om=f(P)\d P_1\wedge ... \wedge \d P_{n-1}.
\]
The representation (\ref{poincare}) 
leads to a hodograph transformation between solutions 
to (\ref{laplace}) and solutions to a system of 
$n$ first order PDEs:
Define $n$ holomorphic $n$--forms in an open set of $\C^{2n}$ by
\be
\label{diff_form}
\omega_{i}=\d x_i\wedge(\ast\d V-\om). 
\ee
The Laplace equation with $|\d V|\neq 0$ is equivalent to $\omega_i=0$.
Selecting an $n$-dimensional surface (an integral manifold) 
in $\C^{2n}$ with $x_i$ as the
local coordinates, and eliminating $P_a$ by cross-differentiating
would lead back to (\ref{laplace}). We are however free to make another
choice and use 
$(V, P_a)$ as local coordinates. This yields (\ref{new_equ}).

Conversely, if $x_i(V, P_a)$ satisfy (\ref{new_equ}) then transforming
it back to (\ref{diff_form}), and closing the ideal 
we deduce that $V(x_i)$ is harmonic.
\koniec
In the three--dimensional case 
$(\Sm, \om)$ is a holomorphic symplectic manifold with local coordinates 
$P_a=(P, Q)$.
We have therefore given an alternative proof
of the equivalence between solutions to the Laplace equation and
the SDiff$(\Sm)$ Nahm equation.
\begin{col}{\em{\cite{W90}}}
\label{sdn}
In the three dimensions the Laplace equation is generically equivalent to the
{\em SDiff$(\Sm)$} Nahm equations
\be
\label{nahmIN}
\dot{x}_1 = \{x_2, x_3\}\qquad
\dot{x}_2 = \{x_3, x_1\}\qquad
\dot{x}_3 = \{x_1, x_2\},
\ee
where $\{, \}$ is a Poisson structure determined by symplectic form $\om$.
\end{col}
The above result fits into a general scheme of 
integrable background geometries \cite{C01}. In our case 
the background geometry is flat,
and this makes the hodograph transformation $x_i(P,Q,V)\rightarrow V(x_i)$ so 
effective.

From now on we shall restrict to the case $n=3$.
Equations (\ref{nahmIN}) are invariant under an infinite--dimensional
group of holomorphic symplectomorphisms of $\Sm$ acting on the leaves of the
foliation $\C\times\Sm$ (the Lax formulation of these equations is given in Appendix B). 
The next result
characterises the solutions for which this symmetry reduces to 
$SL(2, \C)$.
\begin{prop}
\label{CCONE}
Let $x_i(P, Q, V)$ be solutions to the Nahm's system {\em (\ref{nahmIN})}  
with the gauge group 
$SL(2, \C)\subset$ {\em SDiff}$(\Sm)$. Then $V=V(x_i)$ is a harmonic function constant on a central quadric. 
\end{prop}
{\bf Proof.}
Consider a symplectic action of $SL(2, \C)$ on $\Sm$, generated by 
Hamiltonian vector fields $L_i$, such that the Lie brackets satisfy
$[L_i, L_j]=-(1/2)\varepsilon_{ijk}L_k$. Let $h_i=h_i(P, Q)$ be the
corresponding Hamiltonians which satisfy
\[
\{h_i, h_j\}=\frac{1}{2}\varepsilon_{ijk}h_k.
\]
We notice that
\[
\{h_1^2+h_2^2+ h_3^2, h_i\}=0\qquad i=1, 2, 3.
\]
The Poisson structure $\{, \}$ comes from a symplectic two-form,
and so is non-degenerate. Therefore  the Hamiltonians satisfy 
the algebraic constraint
\be
\label{aconstraint2}
h_1^2+h_2^2+h_3^2=C,
\ee
where $C$ is a constant (which can always be scaled to $0$ or $1$). 
The case $C=0$ corresponds to a linear
action, and $C\neq 0$ to a M\"obius action.

Now assume that 
$
x_i=x_i(V, P, Q),
$
satisfy $SL(2, \C)$ Nahm's equation. Therefore 
$x_i=A_{ij}(V)h_j(P, Q)$. The matrix $A_{ij}$ can  be 
made diagonal and so \[
x_1=w_1(V)h_1(P, Q), \qquad
x_2=w_2(V)h_2(P, Q), \qquad x_3=w_3(V)h_3(P, Q).
\] 
The constraint
(\ref{aconstraint2}) implies
\[
\frac{{x_1}^2}{{w_1(V)}^2}+\frac{{x_2}^2}{{w_2(V)}^2}
+\frac{{x_3}^2}{{w_3(V)}^2}=C,
\]
and the level sets of $V$ are quadrics.
\koniec
We conclude that for solutions constant on quadrics
the Nahm's equations (\ref{nahmIN}) reduce to 
the Euler equations.
\be
\label{eulerEQ}
\dot{w}_1=w_2w_3,\qquad \dot{w}_2=w_1w_3,\qquad
\dot{w}_3=w_1w_2.
\ee
These equations readily reduce to
\be
\label{elliptic}
(\dot{w}_3)^2=({w_3}^2+A)({w_3}^2+B),
\ee
where $A={w_1}^2-{w_3}^2$ and $B={w_2}^2-{w_3}^2$ are constants.
Setting 
\[
H(V)= {w_3}^2(V)+\beta_3,
\]
where $\beta_3=(A+B)/3$ yields
\[
\dot{H}^2=(H-\beta_1)(H-\beta_2)(H-\beta_3).
\]
which is (\ref{triangle}) with $n=3$, and 
$\beta_1=\beta_3-A, \beta_2=\beta_3-B$.
\section{Twistor description}
Let  $(x_1, x_2, x_3)$ be coordinates on $\C^3$ in which
$\eta=\d {x_1}^2+\d {x_2}^2+\d {x_3}^2$.
The twistor space ${\cal Z}$ of $\C^3$ is the space 
of all planes $Z\subset \C^3$ with are null with respect to
$\eta$ \cite{J85,W89}. It is the two-dimensional complex manifold 
${\cal Z}=T\CP^1$, which can be seen as follows: null vectors
in $\C^3$ can be parametrised by $k=(1-\ll^2, \sqrt{-1}(1+\ll^2), 2\ll)$, where
$\ll\in\CP^1$. Points of ${\cal Z}$ correspond to 
null 2-planes in $\C^3$ via the incidence relation
\be
\label{section}
\mu=(x_1+\sqrt{-1} x_2)+2x_3\ll-(x_1-\sqrt{-1} x_2)\ll^2.
\ee
Fixing $(\mu , \ll)$
defines a null plane in $\C^3$ with $k$ as its normal. 
An alternate interpretation of (\ref{section}) is to fix $x_i$.
This determines $\mu$ as a function of $\ll$ i.e.\ a section $L_x$ of
${\cal O}(2)\rightarrow \CP^1$.
 
The total space of $T\CP^1$ is equivalent to the total space of the
line bundle $\OO(2)$, and so every holomorphic section of 
${\cal Z}\rightarrow \CP^1$ can
be written as a polynomial quadratic in $\ll$ with complex
coefficients.
The outlined twistor correspondence can be summarised as follows
\begin{eqnarray*}
&&\mbox{points}\longleftarrow\longrightarrow \mbox{holomorphic sections}\\
&&\mbox{null planes}\longleftarrow\longrightarrow\mbox{points}.
\end{eqnarray*}
Another way of defining ${\cal Z}$ is by the double fibration
\be
\label{doublefib}
{\C^3}\stackrel{p}\longleftarrow 
{\cal F}\stackrel{q}\longrightarrow {\cal Z}.
\ee
The correspondence space ${\cal F}=\C^3\times\CP^1$ has a natural fibration
over $\C^3$, and the projection $q:{\cal F}\longrightarrow {\cal Z}$ is
a quotient of ${\cal F}$ by 
the two-dimensional distribution of vectors tangent to a null plane.
In concrete terms this distribution is spanned by
\[
L_0=(\p_1+\sqrt{-1}\p_2)-\ll\p_3, \qquad L_1=\p_3-\ll(\p_1-\sqrt{-1}\p_2).
\]
This leads to an alternative definition
${\cal F}=\{(x,Z)\in\C^3\times {\cal Z}|{Z\in L_x}\}$.

The integral formula for solutions to the Laplace equation
\[
\frac{\p^2 V}{\p {x_1}^2}+\frac{\p^2 V}{\p {x_2}^2}
+\frac{\p^2 V}{\p {x_3}^2}=0
\]
may be elegantly expressed in the twistor terms. Given a holomorphic function
$F(\ll, \mu)$ on ${\cal Z}$ (secretly an element of 
$H^1(\CP^1, \OO(-2))$) restrict it to a section (\ref{section}),
and pull it back to ${\cal F}$.
The general
harmonic function in $\C^3$ is then given by
\be
\label{tf}
V(x)=\oint_{\Gamma} q^*(F(\ll, \mu(\ll)))\d \ll,
\ee
where $\Gamma\subset L_x\cong\CP^1$ is a real closed contour.

It is natural to ask for a characterisation of twistor functions which
give rise to harmonic functions constant on a central quadric. 
We shall first note
that the dilation vector field $\Upsilon={\bf r}\cdot\nabla$ corresponds to 
a holomorphic vector field $\mu/\p \mu$ on ${\cal Z}$. 
Therefore if a twistor function 
$F$ gives rise  to $V$ constant on a central quadric, then $\hat{V}$ 
given by (\ref{exp_V}) with $n=3$ can be written as
\be
\label{tfexp}
\hat{V}(x)=\oint_{\Gamma} q^*\Big(\mu\frac{\p F(\ll, \mu(\ll))}
{\p \mu}\Big)\d \ll.
\ee

A direct attempt to invert (\ref{tf}) with $V$ as in (\ref{exp_V}) leads to a 
messy calculation with an inconclusive outcome.
We shall therefore choose a different route based on 
holomorphic line bundles over ${\cal Z}=\OO(2)$ with the vanishing first Chern
class.
These objects are classified
by elements of $H^1(\OO(2), \OO)$. Let $L$ be such a line bundle corresponding
to  a patching function $\exp{(f)}$, where $f\in H^1(\OO(2), \OO)$. One can view $L$ in two different ways
\begin{enumerate}
\item A pull-back of the cohomology class $\p f/\p\mu$ 
to the correspondence space
gives a solution to the Laplace equation.
\item Lifts of holomorphic sections of $\OO(2)\longrightarrow \CP^1$ to $L$
are rational curves with normal bundle $\OO(1)\oplus\OO(1)$. Moreover
$L$ is fibred over $\CP^1$, and the canonical bundle of $L$ is 
isomorphic to a pull-back of $\OO(-4)$ form $\CP^1$. Therefore $L$ is 
a twistor space of a (complexified) 
hyper-K\"ahler four-manifold $M$ \cite{Pe76, AHS78}.
\end{enumerate}
\begin{theo}
\label{Twistor_th}
Let $f\in H^1(\OO(2), \OO)$ define a holomorphic line bundle
$L\longrightarrow\OO(2)$ with $c_1(L)=0$, and let $\mu=\mu(\ll)$ be a section of
$\OO(2)$ given by {\em(\ref{section})}. The following statements are equivalent
\begin{enumerate}
\item There exists a  
homomorphism of vector bundles,
\be
\label{homomorphism}
\alpha: L \otimes {\bf sl}(2, \C)\longrightarrow TL
\ee
such that rank$(\alpha)=2$. 
\item 
Let $F\in H^1(\OO(2), \OO(-2))$ satisfy 
\[
\mu\frac{\p F}{\p \mu}=\frac{\p f}{\p \mu}.
\]
Then
\be
\label{lap_proof}
V=\oint_{\Gamma} q^*(F(\ll, \mu))\d \ll
\ee
is a solution to the Laplace equation constant 
on a central quadric. 
\end{enumerate}
\end{theo}
The proof of this result is based on a fact that a hyper-K\"ahler 
metric corresponding to the twistor space $L$ can be explicitly given
in terms of a harmonic function $\hat{V}$. We shall postpone this
proof to the next section, where we have characterised 
hyper-K\"ahler metrics corresponding to (Euler derivatives of) harmonic functions  
constant on quadrics.
\section{$SL(2, \C)$ invariant Gibbons-Hawking metrics}
In this section we shall show that harmonic functions constant on 
central quadrics (acted on by the homogeneity operator)
characterise $\C^*$ invariant hyper-K\"ahler metrics 
\cite{GH78} which belong to the BGPP class \cite{BGPP78}.
The construction is a consequence of a simple observation:
If $h$ is a left-invariant metrics on $SL(2,\C)$, and  $\C^*\subset
SL(2, \C)$ then  the metric induced on $\Sm=SL(2, \C)/\C^*$ is conformal
to a metric on a complex quadric.

Recall that a four-dimensional manifold $M$ is 
(complexified) hyper-K\"ahler if it admits K\"ahler 
structures $\Omega_i, i=1, 2, 3$ compatible with a fixed (holomorphic) 
Riemannian
metric $g$ and such that the endomorphisms $I_i$ given by 
$g(I_iX,Y)=\Omega_i(X,Y)$
satisfy the algebraic relation of quaternions. 
The hyper-K\"ahler  metrics  on
$\C\times SL(2,\C)$
with a transitive tri-holomorphic action
of $SL(2, \C)$  can be put in the form \cite{BGPP78}
\be
\label{BGPP}
g=w_1w_2w_3\d V^2+\frac{w_2w_3}{w_1}(\sigma_1)^2
+\frac{w_1w_3}{w_2}(\sigma_2)^2+\frac{w_1w_2}{w_3}(\sigma_3)^2,
\ee
where $w_1, w_2, w_3$ are functions of $V$, and
$\sm_i$ are  left invariant one-forms on $SL(2, \C)$ which satisfy 
\be
\label{cijk}
\d \sigma_1=\sigma_2\wedge\sigma_3,\qquad
\d \sigma_2=\sigma_3\wedge\sigma_1,\qquad
\d\sigma_3=\sigma_1\wedge\sigma_2.
\ee
The self-dual two-forms are
\begin{eqnarray}
\label{BGPPom}
\Om_1&=&w_1\sigma_2\wedge\sigma_3+w_2w_3\sigma_1\wedge\d V,\nonumber\\
\Om_2&=&w_2\sigma_3\wedge\sigma_1+w_1w_3\sigma_2\wedge\d V,\nonumber\\
\Om_3&=&w_3\sigma_1\wedge\sigma_2+w_1w_2\sigma_3\wedge\d V,
\end{eqnarray}
and the hyper-K{\"a}hler condition $\d \Om_i=0$ 
is\footnote{The $SL(2,\C)$ action fixes all complex structures, 
so the invariant frame is covariantly constant.} equivalent to the Euler
equations (\ref{eulerEQ}).
Therefore $A={w_1}^2-{w_3}^2$ and $B={w_2}^2-{w_3}^2$ are constants.
The BGPP metric (\ref{BGPP}) is flat if $A=B=0$, and is never complete
if $AB(A-B)\neq 0$. The remaining cases correspond to   
a complete metric known as the Eguchi--Hanson solution.

Now  choose a one-dimensional 
subgroup  $\C^*\subset SL(2, \C)$.
Hyper-K{\"a}hler metric with a tri-holomorphic $\C^*$ action can
be put in the Gibbons--Hawking  form \cite{GH78}
\be
\label{GH}
g=\hat{V}(\d {x_1}^2+\d {x_2}^2+\d {x_3}^2)+\hat{V}^{-1}(\d T+A)^2,
\ee
where $\ast \d \hat{V}=\d A$. Here $K=\p/\p T$ generates the $\C^*$ action,
and $x_i$ are defined up to addition of a constant
by $\d x_i=K\hook \Om_i$. The one-form $A$ is defined on the orbits of $K$.

To characterise the harmonic functions $\hat{V}$ for which (\ref{GH})
belongs to the BGPP class we could  expand $K$
in a left-invariant basis of $SL(2, \C)$, eliminate
the Euler angles and $V$  in favour of $(x_i, t)$ and observe that
$\hat{V}^{-1}=g(K, K)$, where $g$ is given by (\ref{BGPP}).
This is essentially done in \cite{GOVR88}  and more explicitly in \cite{G03}, 
where the $n=3$ case of (\ref{exp_V}) is obtained. 

We shall adopt a variation of this  approach,
and use the Nahm equations to establish the following

\begin{prop} The Gibbons--Hawking metric {\em (\ref{GH})}
belongs to the BGPP class {\em(\ref{BGPP})} 
if and only if $\hat{V}={\bf r}\cdot\nabla V$, and $V$ is a harmonic function
constant on a central quadric.
\label{Gibbons_Hawking}
\end{prop}
{\bf Proof.} Let $(M=\C^*\times SL(2, \C), g)$ be 
a hyper-K\"ahler four-manifold with a transitive
and tri-holomorphic $SL(2, \C)$ action, and
let $\gamma:SL(2, \C)\longrightarrow \Sm=SL(2, \C)/\C^*$ be a 
complexified Hopf bundle.

The Corollary \ref{sdn}  allows us to 
introduce $(P, Q, \hat{V}, T)$ as local coordinates on $M$ 
such that $(P, Q)$ are local coordinates on $\Sm$, 
$\hat{V}$ is a coordinate on the fibres
of $\gamma$, and $T$ parametrises the $SL(2, \C)$ orbits. 
We shall regard $x_i=x_i(P, Q, \hat{V})$ as  functions on $\C^3$.
The hyper--K\"ahler condition is  then equivalent to SDiff$(\Sm)$ Nahm
equation 
\[
\frac{\p x_i}{\p \hat{V}}=\frac{1}{2}\varepsilon_{ijk}\{x_j, x_k\},
\]
where the Poisson brackets are taken 
with respect to the symplectic form $\om=\d A$.

Consider the left action of $SL(2, \C)$ on itself, generated 
by left invariant vector fields $L_i$, such that
\[
{\cal L}_{L_i}g=0,\qquad  [L_i, L_j]=-\frac{1}{2}\varepsilon_{ijk}L_k,\qquad 
L_i\hook\sigma_j=\delta_{ij}.
\]
The push-forward vector fields $\gamma_{\ast}(L_i)$
generate symplectomorphisms of $\Sm$, and so they correspond to 
the Hamiltonians $h_i=h_i(P,Q)$ (these are the Hamiltonians used in the proof
of Proposition \ref{CCONE}), such that
\be
\label{left_right}
\gamma_\ast(L_i)(h_j)=\{h_i, h_j\}=\frac{1}{2}\varepsilon_{ijk}h_k,
\ee
and the algebraic constraint (\ref{aconstraint2}) holds. 
Let $K=\gamma^{\ast}(h_1)L_1+\gamma^{\ast}(h_2)L_2+\gamma^{\ast}(h_3)L_3 $. Note that
\[
[K, L_i]=0
\]
as a consequence of (\ref{left_right}). 
The moment-maps $x_i$ are given by
\[
\d x_i=
K\hook\Om_i=\gamma^{\ast}(h_i)\dot{w}_i\d V+w_i\d (\gamma^{\ast}(h_i))
=\gamma^{\ast}\d(w_i h_i), \qquad\mbox{(no summation)},
\]
where $\Om_i$ are given by (\ref{BGPPom}).
The formula (\ref{GH}) implies that $\hat{V}=(g(K, K))^{-1}$.
Making a substitution $w_i=\sqrt{H-\beta_i}$,
where $H(V)$ satisfies  (\ref{triangle}) with $n=3$ reveals that $\hat{V}$
is given by formulae (\ref{exp_V}) with $n=3$. Therefore $\hat{V}=
\Upsilon(V)$ is an Euler derivative of a solution 
(\ref{new_sol}) constant on a central 
quadric.
\koniec
We are now ready to present a proof of Theorem \ref{Twistor_th}, and give 
the characterisation of inverse twistor 
functions corresponding to $V$ constant on a central quadric. 

{\bf Proof of Theorem \ref{Twistor_th}.}
From (\ref{dimh1}) it follows that cocycles in $H^1(\CP^1, \OO(n))$ can be
represented by coboundaries if $n\geqslant -1$. The freedom one has is measured
by $H^0(\CP^1, \OO(n))$. In particular 
$H^1(\CP^1, \OO(-1))=0$ and its cocycles can be uniquely represented as 
coboundaries. Let $\pi=(\pi_0, \pi_1)$ be homogeneous coordinates on $\CP^1$.
Let  $U_0$ and  $U_1$ be a covering of ${\OO(2)}$ such that 
$\pi_{1}\neq 0$ on $U_0$, and $\pi_{0}\neq 0$ on 
$U_1$, and
let $q^*U_{\alpha}$ denote the open sets
on the correspondence space $\C^3\times\CP^1$ that are the pre-image of
$U_{\alpha}$, for ${\alpha}=0,1$.
Let $f\in H^1(\OO(2), \OO)$ be a logarithm of a patching function of $L$. Consider
\[
\pi\frac{\p f}{\p \mu}\in \C^2\otimes H^1(\OO(2), \OO(-1)),
\]
which is  homogeneous of degree $-1$.
Then we  restrict it to a
section of $\OO(2)\longrightarrow \CP^1$ and pull it back to ${\cal F}$
by (\ref{doublefib}),
where we can split as $\pi\p f/\p\mu=h_0-h_1$.

Here $h_{\alpha}$ is  holomorphic on 
$q^*{U_\alpha}$, and is given by
\be
\label{sparf}
h_{\alpha}
=\oint_{\Gamma_\alpha}\frac{\rho}{\pi\cdot\rho} 
\frac{\p f}{\p \mu}
\rho \cdot \d\rho,
\ee
where $\rho$ are homogeneous coordinates on $\CP^1$, and 
$\rho\cdot\d\rho=\d \ll$ in affine coordinates.
The contours $\Gamma_\alpha$ and 
are homologous to the equator of $\CP^1$ in $U_0\cap  U_1$ and
are such that $\G_0- \G_1$ surrounds the point
$\rho=\pi$. We see that $\pi\p f/\p\mu=h_0-h_1$ 
follows from the Cauchy's integral formula. Moreover $\pi\cdot(\pi f)=\pi\cdot h_0-\pi\cdot h_1 =0$.
Therefore $\hat{V}:=\pi\cdot h_0=\pi\cdot h_1$ is a global holomorphic function
homogeneous of degree $0$, and so by the Liouville Theorem 
it is constant on $\CP^1$. 
The formula (\ref{sparf}) implies that $\hat{V}$ is a harmonic
function  explicitly given by 
\be
\label{lap_proof1}
\hat{V}=\oint_{\Gamma} q^*\Big(\frac{\p f(\ll, \mu)}{\p \mu})\Big)\d \ll,
\ee
which is (\ref{tfexp}) with 
\[
\mu\frac{\p F}{\p \mu}=\frac{\p f}{\p \mu}.
\]
Now we shall regard $L$ as a twistor space of  a hyper-K\"ahler manifold
$(M, g)$. The fibres of $L\longrightarrow \CP^1$ are symplectic manifolds, 
where the symplectic two-form $\Om$ takes values in $\OO(2)$.
The Hamiltonian vector field of $\mu\in\Gamma(\OO(2))$ with respect to $\Omega$ 
preserves the projection of $L$ onto $\CP^1$, and therefore
induces  a tri-holomorphic Killing vector $K=\p/\p T$ on $M$. 
Metrics which admit such Killing vectors are locally given by 
(\ref{GH})\cite{TW79}. 
We shall show  that $\hat{V}$ appearing in (\ref{GH}) is indeed
the same as (\ref{lap_proof1}).

 Introduce local homogeneous coordinates
$(\pi, \mu, \zeta_\alpha)$
on each set $U_\alpha$ of some Stein cover of the twistor space; here 
$\zeta_\alpha$ is a
fibre coordinate up the fibres of the affine line bundle $L\rightarrow \OO(2)$ 
on $U_\alpha$ with patching relations
$\zeta_0=\zeta_1+f$ on $U_0\cap U_1$.  In these coordinates
$$
\Om=\d_h \mu\wedge\d_h\zeta_\alpha\, ,$$ where $\d_h$ denotes the exterior
derivative in which $\pi$ is held constant. 
The two--form
$\Om$ is globally defined on vector fields tangent to the fibres of
$L\longrightarrow \CP^1$ as $f$ does not depend on
$\zeta_i$.
In order to calculate the self-dual  two-forms $\Om_i$ (and so the metric)  
associated to  $f$ we
pullback $\Om$ to 
$M\times\CP^1$, and determine $\zeta_{\alpha}$ using 
integral splitting formulae analogous to  (\ref{sparf})
(note however, that $h_\alpha$ are not twistor functions, as they don't 
descend from ${\cal F}$ to ${\cal Z}$. On the other  hand $\zeta_\alpha$
are twistor function for $L$). 
This yields
\begin{eqnarray*}
q^*(\Om)&=&(\Om_1+\sqrt{-1} \Om_2)+2\Om_3\ll-(\Om_1-\sqrt{-1} \Om_2)\ll^2,\qquad
\mbox{where}\\ 
\Om_i&=&(\d T+A)\wedge\d x_i-\frac{1}{2}{\hat{V}}\varepsilon_{ijk}\d x_j\wedge\d x_k,
\end{eqnarray*}
The two-forms $\Om_i$ are the self-dual two forms of (\ref{GH}), and 
$x_i$ are given by (\ref{section}).
The function $\hat{V}$, and a one-form  $A=A_i\d x^i$ are 
given in terms of a $2\times 2$ matrix
\[
\Phi_{BC}=
\oint_{\Gamma_\alpha}\frac{\rho_B\iota_C}{\rho\cdot \iota} 
\frac{\p f}{\p \mu}
\rho \cdot \d\rho,\qquad \iota=(\iota_0, \iota_1)\in\CP^1, \quad B, C=0, 1
\]
by
\[
\Phi=
\left (
\begin{array}{cc}
A_1+\sqrt{-1}A_2&A_3+\hat{V}\\
A_3-\hat{V}&-(A_1-\sqrt{-1}A_2)
\end{array}
\right ).
\]
Therefore $\hat{V}$ is as in  (\ref{lap_proof1}), and the monopole equation
$\ast\d \hat{V}=\d A$ is automatically satisfied.
The two-forms $\Om_i$ are closed, as a consequence of the monopole equation.

Proposition \ref{Gibbons_Hawking}
implies that  $\hat{V}={\bf r}\cdot\nabla V$, where $V$ 
is harmonic and constant on a central quadric iff the associated metric  
(\ref{GH}) belongs to the BGPP class.
To complete the proof we need to characterise the structures on $L$
induced by a tri-holomorphic $SL(2, \C)$ action on $M$.
The $SL(2, \C)$ action on $M$  preserves $g$ therefore it induces a 
holomorphic $SL(2, \C)$ action on $L$. Each Killing vector $L_i$ on $M$
induces a 
holomorphic vector field $\hat{L}_i$ on $L$, and this gives rise \cite{H95}
to a homomorphism
(\ref{homomorphism}) such that $2\leq\,$rk$\,(\alpha)\leq 3$. The
group action on $M$ is tri-holomorphic,
therefore $Q_i=\hat{L}_i\hook\d \ll=0$ (in general $Q_i\in \OO(2)$ gives rise 
to a vector field Hamiltonian with respect to $\d \pi_0\wedge\d\pi_1$ 
which rotates
the self-dual two-forms).
The orbits of the $SL(2, \C)$ action on 
$L$ are contained in the two-dimensional fibres of $L\longrightarrow \CP^1$.
We conclude that rank$(\alpha)=2$.
\koniec 
\section{Example}
We shall illustrate Proposition (\ref{Gibbons_Hawking}) with the example
of a harmonic function corresponding  to the Eguchi--Hanson metric.  
We impose the Euclidean reality conditions, and work
with real harmonic functions on $\R^3$ giving rise to real
hyper-K\"ahler metrics. 

The one-forms $\sigma_i$ can be explicitly given in terms of Euler
angles
\[
\sigma_1=\cos{\psi}\;\d \theta+\sin{\psi}\sin{\theta}\;\d\phi,\qquad
\sigma_2=\sin{\psi}\;\d \theta-\cos{\psi}\sin{\theta}\;\d\phi,\qquad
\sigma_3=\d\psi+\cos{\theta}\;\d \phi.
\]
The functions
\[
h_1=\sin{\theta}\sin{\psi},\qquad h_2=-\sin{\theta}\cos{\psi},\qquad
h_3=\cos{\theta}
\]
satisfy (\ref{left_right}) with
the symplectic form $\om=\d(\cos{\theta})\wedge\d\psi$.  

Consider the
$SU(2)$ invariant hyper--K\"ahler metric (\ref{BGPP}) with
$w_1=w_2\neq w_3$. The Euler equations (\ref{eulerEQ}) yield
$w_3=\rho(V), w_1=w_2=\sqrt{\rho^2-a^2}$, and (using $\rho$ as a
coordinate)
\be
\label{EHmetric}
g=\frac{\rho}{\rho^2-a^2}
\d \rho^2+\rho(\sm_1^2+\sm_2^2)+\frac{\rho^2-a^2}{\rho}\sm_3^2.
\ee
This is the metric of Eguchi and Hanson \cite{EH78}. It is complete,
as the apparent singularity at $\rho=a$ is removed by allowing
\[
a^2<\rho^2,\qquad 
 0\leq\psi\leq 2\pi,\qquad 0\leq \phi\leq 2\pi,\qquad 0\leq\theta\leq\pi.
\]
Identifying  $\psi$ modulo $2\pi$ makes the  surfaces 
of constant $\rho^2>a^2$ into $\RP^3$.  At large value of $\rho^2$ the metric is
asymptotically locally Euclidean.

The metric  (\ref{EHmetric}) can be put in the Gibbons--Hawking form with
$\hat{V}=\rho/(\rho^2-a^2\cos^2{\theta})$. 
To see this perform a  coordinate transformation $x_i=\om_ih_i,
T=\phi+\psi$ which yields (\ref{GH}) with 
\[
\hat{V}=|{\bf r+a}|^{-1}+|{\bf r-a}|^{-1},\qquad
\mbox{where}\qquad {\bf a}=(0, 0, a).
\]
We verify that ${\bf r}\cdot\nabla V=\hat{V}$,
where the harmonic function
\[
V=-\frac{2}{a}\mbox{arccoth}{\Big(\frac{|{\bf r+a}|+|{\bf
r-a}|}{2a}\Big)}
\]
is obtained from (\ref{new_sol}) with $n=3$, and
$\beta_1=\beta_2=a^2, \beta_3=0$. 
The potential $V$  is constant on the ellipsoid
\[
\frac{x^2+y^2}{a^2((\coth{(aV/2)})^2-1)}+\frac{z^2}{a^2(\coth{(aV/2)})^2}=1,
\]
in agreement with Proposition (\ref{Gibbons_Hawking}).
\section*{Acknowledgements}
I thank Gary Gibbons and Paul Tod for useful discussions.
This research was partly supported by NATO grant PST.CLG.978984.
\section*{Appendix A. Bundles over $\CP^1$}
\appendix
\def\theequation{\thesection{A}\arabic{equation}}
Let $\C^2$ be  a symplectic vector space, with
anti-symmetric product
\[
\pi\cdot\rho=\pi_0\rho_1-\pi_1\rho_0=-\rho\cdot\pi,  
\]
where $\pi=(\pi_0, \pi_1), \rho=(\rho_0, \rho_1)\in \C^2$.
Remove
$\pi=(0, 0)$ and use $\pi$ as homogeneous coordinates on $\CP^1$.
We shall also use the affine coordinate $\ll=\pi_{0}/\pi_{1}$.
Holomorphic functions on $\C^2-0$ extend to holomorphic functions on $\C^2$
(Hartog's Theorem). Therefore homogeneous functions on $\CP^1$ are polynomials.
In particular, holomorphic functions homogeneous of degree $0$ 
are constant (Liouville theorem).
Let us summarize some facts about holomolrphic line bundles
over $\CP^1$. First define a tautological line bundle 
\[
\OO(-1)=\{(\ll, (\pi_0, \pi_1))\in \CP^1\times\C^2|\ll=\pi_0/\pi_1\}.
\]
Other line bundles can be obtained from $\OO(-1)$ by algebraic operations:
\[
\OO(-n)=\OO(-1)^{\otimes n}, \qquad \OO(n)=\OO(-n)^*, \qquad\OO=\OO(-1)\otimes\OO(1), 
\qquad n\in \NN.
\]
Equivalently  ${\cal O}(n)$ denotes the line bundle over $\CP^1$ with transition functions
$\ll^{-n}$ from the set $\ll\neq\infty$ to $\ll\neq 0$ (i.e.\ Chern class
$n$).  Its sections
are given by functions homogeneous of degree $n$ in a sense that
$
f(\xi\pi)=\xi^nf(\pi).
$
These are polynomials in $\ll$ of degree $n$ with complex 
coefficients.
The theorem of Grothendick states that all holmorphic line bundles
over a rational curve are equivalent to $\OO(n)$ for some $n$.
The spaces of global  sections, and the first cohomology groups are
\be
 H^0(\CP^1,{\cal O}(n))=\left\{ \begin{array}{ll}
                       0    & \mbox{ for}\; n<0 \\
                       \C^{n+1} & \mbox{for}\; n\geq 0.
                       \end{array}
                           \right.\qquad
\label{dimh1}
 H^1(\CP^1,{\cal O}(-n))=\left\{ \begin{array}{ll}
                       0    & \mbox{ for}\; n<2 \\
                       \C^{n-1} & \mbox{ for}\; n\geq 2 .
                       \end{array}
                           \right.
\ee
\section*{Appendix B. SDiff$(\Sm)$ Nahm's equations}
\appendix
\def\theequation{\thesection{B}\arabic{equation}}
Let $G$ be a Lie group  
and let  $[\;,\;]$ be the Lie bracket in the corresponding Lie algebra 
${\bf g}$, 
The Nahm equations for three ${\bf g}$-valued functions $X_i=X_i(V)$ are
\be
\label{nahm}
\dot{ X_1} = [X_2,X_3],\qquad
\dot{ X_2} = [X_3,X_1],\qquad
\dot{ X_3} = [X_1,X_2].
\ee
These equation admit a Lax representation.
Let
\[
A(\ll)=(X_1+\sqrt{-1} X_2)+2X_3\ll-(X_1-\sqrt{-1} X_2)\ll^2.
\]
Then
\begin{eqnarray}
\label{lax_nahm}
\dot{A} &=& [X_2-\sqrt{-1} X_1,X_3] +2[X_1,X_2]\ll-[X_2+\sqrt{-1} X_1,X_3]\ll^2\nonumber \\
 &=& [A,-\sqrt{-1} X_3+\sqrt{-1} (X_1-\sqrt{-1} X_2)\ll]\nonumber \\
 &=& [A, B],
\end{eqnarray}
where $B=-2\sqrt{-1}X_3+\sqrt{-1}(X_1-\sqrt{-1}X_2)\ll$. 
The matrix $A(\ll)$
should be thought of as an $\OO(2)$ valued section of a two-dimensional
complex vector bundle over $\CP^1$. 
Let $(\mu, \ll)$ be the local coordinates on the total space
of $\OO(2)$. 
The zero locus of the  characteristic equation
\[
S=\{(\mu, \ll)\in{\cal Z}| \det({\bf 1}\mu-A(\ll))=0\}
\]
defines an algebraic curve $S$, called the spectral curve of $A(\ll)$ 
which (as a consequence of (\ref{lax_nahm})) is preserved by
the Nahm equations. In the $SL(2, \C)$ case  the Riemann surface $S$ has genus one, and is 
a torus parametrised by the elliptic function (\ref{elliptic}).

For each point $({\mu, \ll})$ on $S$ we have a one-dimensional
space $L_{({\mu, \ll})}=\ker \,({\bf 1}\mu-A(\ll))$
and this varies with $V$. It forms a line bundle over the
spectral curve. Hitchin \cite{H83} shows that this line bundle evolves along a straight line on the Jacobian.

Now assume that ${\bf g}$ is the infinite-dimensional Lie algebra 
sdiff$(\Sm)$ of holomorphic symplectomorphisms of a two-dimensional complex 
symplectic manifold $\Sm$ with local holomorphic 
coordinates $P, Q$ and the holomorphic symplectic
structure $\om$. Elements of
sdiff$(\Sm)$ are represented by the Hamiltonian vector fields $X_h$
such that
$
X_h\hook \om=\d h
$ 
where $H$ is a $\C$-valued function on $\Sm$.
The Poisson algebra of functions which we are going to use 
is homomorphic to sdiff$(\Sm)$. 
We shall make the replacement 
$
[\;,\;]\longrightarrow \{\;,\;\}. 
$
in formulae (\ref{nahm},\ref{lax_nahm}).
Here $\{, \}$ is a Poisson structure defined by $\om$.
The components of
$X_i$ are therefore replaced  by  Hamiltonians $x_i(P, Q, V)$ 
generating the 
symplectomorphisms of $\Sm$.
The Lax representation for the SDiff$(\Sm)$ Nahm system is
\be
\frac{\p\Psi}{\p V}=\{({\mu}/{\ll})_+, \Psi\},\qquad
\{\mu, \Psi\}=0,
\ee
where $\Psi=\Psi (P, Q, V, \ll)$, and
\[
\mu(\ll)=(x_1+\sqrt{-1} x_2)+2x_3\ll-(x_1-\sqrt{-1} x_2)\ll^2,\;\; 
({\mu}/{\ll})_+ =-2\sqrt{-1}x_3+\sqrt{-1}(x_1-\sqrt{-1}x_2)\ll. 
\]
The compatibility conditions for this over-determined system 
yield (\ref{nahmIN}).

\end{document}